\documentclass[a4paper,11pt]{article}
\usepackage{amsmath,enumerate}
\usepackage{amssymb,geometry,latexsym,theorem}
\usepackage[all]{xy}
\usepackage[ps2pdf, bookmarks=true,
            bookmarksnumbered=true, breaklinks=true,
            pdfstartview=FitH, hyperfigures=false,
            plainpages=false, naturalnames=true,
            colorlinks=true,
            pdfpagelabels]{hyperref}
\usepackage[usenames]{color}
{\theorembodyfont{\rmfamily}\newtheorem{example}{Example}[section]}
{\theorembodyfont{\rmfamily}}
{\theorembodyfont{\rmfamily}}
{\theorembodyfont{\rmfamily}}
\newtheorem{prop}[example]{Proposition}
\newtheorem{thm}[example]{Theorem}
{\theorembodyfont{\rmfamily}\newtheorem{rem}[example]{Remark}}
\newtheorem{cor}[example]{Corollary}

\def\I{{\mathbb I}}

\def\Z{\mathbb{Z}}

\def\ES{\mathbb{S}}

\hypersetup{ pdfauthor   = {Ronald Brown}, pdftitle    = {Topology
and Groupoids}, pdfsubject  = {Mathematics}, pdfkeywords = {},
pdfcreator  = {LaTeX with hyperref package}, pdfproducer = {dvips
+ ps2pdf + thumbpdf} }

\usepackage{makeidx}

\newcommand{\Ob}{\operatorname{Ob}}

\def\2{\mathsf{2}}

\def\Proof{\noindent {\bf{Proof}}\ }
\def\sqbox{\rule{0em}{1ex} \hfill $\Box$}\def\Real{\mathbb{R}}
\def\phi{\varphi}

\begin{document}
\title{Groupoids, the Phragmen-Brouwer Property, \\ and the Jordan
Curve Theorem}
\author{Ronald Brown\\ University of Wales, Bangor\thanks{www.bangor.ac.uk/r.brown.
This work was partially supported by a Leverhulme Emeritus
Fellowship (2002-2004). }} \maketitle

\begin{abstract}
We publicise a proof of the Jordan Curve Theorem which relates it
to the Phragmen-Brouwer Property, and whose proof uses the van
Kampen theorem for the fundamental groupoid on a set of base
points.
\end{abstract}

\section{Introduction}
This article extracts from \cite{Brown15} a proof of the Jordan
Curve Theorem based on the use of groupoids, the van Kampen Theorem
for the fundamental groupoid on a set of base points, and the use of
the Phragmen-Brouwer Property. In the process, we give short proofs
of two results on the Phragmen-Brouwer Property (Propositions
\ref{R:8.5.1}, \ref{R:8.5.3}).  There is a renewed interest in such
classical results\footnote{See also the web site on the Jordan Curve
Theorem: {\url{http://www.maths.ed.ac.uk/~aar/jordan/}}}, as shown
in the article \cite{Siebenmann} which revisits proofs of the
Schoenflies Theorem.

There are many books containing a further discussion of this area.
For more on the Phragmen-Brouwer property, see \cite{Whyburn1} and
\cite{Wilder1}.  Wilder lists five other properties which he shows
for a connected and locally connected metric space are each
equivalent to the PBP.  The  proof we give of the Jordan Curve
Theorem is adapted from \cite{Munkres1}.  Because he does not have
our van Kampen theorem for non-connected spaces, he is forced into
rather special covering space arguments to prove his replacements
for Corollary \ref{cor:retract} (which is originally due to
Eilenberg \cite{Eilenberg}), and for Proposition \ref{R:8.5.1}.

The intention is to make these methods more widely available, since
the  1988 edition of the book \cite{Brown15} has been out of print
for at least ten years, and the new edition is only just now
available.

I mention in the same spirit that the results from \cite{Brown15} on
orbit spaces have been made available in \cite{BrownHigginsorbit}.
As further example of the use of groupoid methods, this time in
combinatorial group theory, is in \cite{Braun}, which gives a new
result combining the Kurosch theorem and a theorem of Higgins which
generalises Grusko's theorem. Coverings of non connected topological
groups are discussed in \cite{BrownMucuk}: essential use is made of
the well known equivalence, for suitable $X$, between the categories
of covering maps over $X$ and of covering morphisms over the
fundamental groupoid $\pi_1 X$. Higgins in \cite{Higgins-graphs}
gives a powerful normal form theorem for what he calls the {\it
fundamental groupoid} of a graph of groups, avoiding the usual
choice of a base point or a tree.

Note that we use groupoids not to give nice proofs of  theorems on
the fundamental {\it group of a space with base point}, but because
we maintain that   theorems in this area  are about the fundamental
{\it groupoid on a set of base points}, where that set is chosen in
a way appropriate to the geometry of the situation at hand. The set
of objects of a groupoid gives a spatial component to group theory
which allows for more powerful and more easily understood modelling
of geometry, and hence for more computational power. Indeed this was
the message of the paper \cite{Brown5} and even the first 1968
edition of \cite{Brown15}.

I would like to thank Michel Zisman for significant improvements to
parts of the exposition, and a referee for helpful comments.

\section{The groupoid van Kampen theorem}
We assume as known the notion of the fundamental groupoid
$\pi_1XJ$ of a topological space $X$ on a set $J$: it consists of
homotopy classes rel end points of paths in $X$ joining points of
$J \cap X$. We say the pair $(X,J)$ is {\it connected} if $J$
meets each path component of $X$. The following theorem was proved
in \cite{Brown5} (see also \cite[6.7.2]{Brown15}).
\begin{thm}[van Kampen Theorem] \label{thm:vKT}
Let the space $X$ be the union of open subsets $U,V$ with
intersection $W$, let $J$ be a set and suppose the pairs $(U,J),
(V,J), (W,J)$ are connected. Then the pair $(X,J)$ is connected
and the following diagram of morphisms induced by inclusion is a
pushout in the category of groupoids: \begin{equation*} \xymatrix{
\pi_1 WJ \ar[r] \ar [d] & \ar [d] \pi_1 VJ \\
\pi_1UJ \ar [r]& \pi_1 XJ }
\end{equation*}
\end{thm}

This has been generalised to  unions of any number of open sets in
\cite{BrownRazakSalleh1}. There then has to be an assumption that
$(U,J)$ is connected for any 3-fold  (and hence also 1- and 2-fold)
intersection $U$ of the sets of the cover.

\section{Pushouts of groupoids}
In order to apply Theorem \ref{thm:vKT}, we need  some
combinatorial groupoid theory. This was set up in \cite{Higgins4},
\cite{Brown15}. We  first explain here how to compute an object
group $H(x)$ of a groupoid $H=G/R$ given as the quotient of a
groupoid $G$ by a totally disconnected graph $R=\{R(x)\mid x \in
\Ob(G)\}$ of relations: of course $G/R$ is defined by the obvious
universal property, and has the same object set as $G$.

Recall from \cite[8.3.3]{Brown15} that:
\begin{prop} \label{prop:quotients} {\rm (a)} If $G$ is a connected groupoid,
and $x \in \Ob(G)$,
then there is  a retraction $r: G \to G(x)$ obtained by choosing
for each $y \in \Ob(G)$ an element $\tau_y \in G(x,y)$, with
$\tau_x=1_x$. \\
{\rm (b)} If further $R=\{R(y)\mid y \in \Ob(G)\}$ is a family of
subsets of the object groups $G(y)$ of $G$, then the object group
$(G/R)(x)$ is isomorphic to the object group $G(x)$ factored by
the relations $r(\rho)$ for all $\rho \in R(y), y \in \Ob(G)$.
\end{prop}

We assume as understood the notion of free groupoid on a
(directed) graph. If $G$, $H$ are groupoids then their free
product $G * H$ is given by the pushout of groupoids
$$\xymatrix{\Ob(G) \cap \Ob(H) \ar [r]^-i \ar [d]_-j  & H \ar [d] \\
G \ar [r] & G*H }$$ where  $\Ob(G) \cap \Ob(H)$ is regarded as the
subgroupoid of identities of both $G,H$ on this object set, and
$i,j$ are the inclusions. We assume, as may be proved from the
results of \cite[Chapter 8]{Brown15}:
\begin{prop}\label{prop:freeproductof free} If $G,H$ are free groupoids, then
so also is $G * H$.
\end{prop}

If $J$ is a set, then by the category of groupoids over $J$ we
mean the category whose objects are groupoids with object set $J$
and whose morphisms are morphisms of groupoids which are the
identity on $J$.

\begin{prop} \label{R:8.4.9}
Suppose given a pushout of  groupoids over $J$
\begin{equation} \label{Eq:8.4:1}
\begin{array}{c}
\xymatrix@=1pc{ C \ar[rr]^*+{i} \ar[dd]_*+{j} &&
A \ar[dd]^*+{u} \\ \\
B \ar[rr]^*+{v} && G }
\end{array} \tag{1}
\end{equation}
such that $C$ is totally disconnected and $A,B$ are connected.  Let
$p$ be a chosen element of $J $. Let $r : A \to A(p)$, $s : B \to
B(p)$ be retractions obtained by choosing elements $\alpha_{x} \in
A(p,x)$, $\beta_{x} \in B(p,x)$, for all $x \in J$, with $\alpha_{p}
= 1$, $\beta_{p} \in 1$.  Let $f_{x} = (u\alpha_{x})^{-1}
(v\beta_{x})$ in $G(p)$, and let $F$ be the free group on the
elements $f_{x}$, $x \in J$, with the relation $f_{p} = 1$.  Then
the object group $G(p)$ is isomorphic to the quotient of the free
product group
$$
A(p) * B(p) * F
$$
by the relations
\begin{equation} \label{Eq:8.4:2}
(r i \gamma)f_{x}(sj\gamma)^{-1}f^{-1}_{x} = 1 \tag{2}
\end{equation}
for all $x\in  J$  and all $\gamma  \in  C(x,x)$.
\end{prop}

\Proof We first remark that the pushout (\ref{Eq:8.4:1}) implies
that the groupoid $G$ may be presented as the quotient of the free
product groupoid $A * B$ by the relations $(i\gamma )(j\gamma
)^{-1}$ for all $\gamma \in C$.  The problem is to interpret this
fact in terms of the object group at $p$ of $G$.

To this end, let $T,S$ be the tree subgroupoids of $A, B$
respectively generated by the elements $\alpha_{x}, \beta_{x}$, $x
\in J$.  The elements $\alpha_{x}, \beta_{x}$, $x \in J$, define
isomorphisms
$$
\phi  : A \to A(p) * T, \qquad \psi  : B \to B(p) * S
$$
where if $g \in  G(x,y)$ then
$$
\phi g  = \alpha_{y}(rg)\alpha^{-1}_{x}, \qquad \psi g
    = \beta_{y}(sg)\beta^{-1}_{x}.
$$
So $G$ is isomorphic to the quotient of the groupoid
$$
H =  A(p) * T * B(p) * S
$$
by the relations
$$
(\phi i\gamma)(\psi j\gamma)^{-1} = 1
$$
for all $\gamma \in C$.  By Proposition \ref{prop:quotients}, the
object group $G(p)$ is isomorphic to the quotient of the group
$H(p)$ by the relations
$$
(r\phi i\gamma)(r\psi j\gamma)^{-1} = 1
$$
for all $\gamma  \in  C$.

Now if $J' = J \setminus \{p\}$, then $T, S$ are free groupoids on
the elements $\alpha_{x}, \beta_{x}$, $x \in J'$, respectively. By
Proposition \ref{prop:freeproductof free}, and as the reader may
readily prove, $T * S$ is the free groupoid on all the elements
$\alpha_{x}, \beta_{x}$, $x \in J'$. It follows from
\cite[8.2.3]{Brown15} (and from Proposition \ref{prop:quotients}),
that $(T
* S)(p)$ is the free group on the elements $r\beta_{x} =
\alpha^{-1}_{x}\beta_{x} = f_{x}$, $x \in J'$.  Let $f_{p}=1 \in
F$. Since
$$
r\phi i\gamma  = ri\gamma, \qquad r\psi j\gamma  = f_{x}(sj\gamma
)f^{-1}_{x},
$$
the result follows. \sqbox
\begin{rem}
The above formula is given in essence in van Kampen's paper
\cite{Kampen1}, since he needed the case of non connected
intersection for applications in algebraic geometry. However his
proof is difficult to follow, and a modern proof for the case of
connected intersection was  given by Crowell in \cite{Crowell1}.
\sqbox
\end{rem}

There is a consequence of the above computation (see
\cite{Eilenberg}) which we shall use in the next section in proving
the Jordan Curve Theorem.

First, if $F$ and $H$ are groups, recall that we say that $F$ is a
{\it retract} of $H$ if there are morphisms $\imath : F \to H$,
$\rho : H \to F$ such that $\rho \imath = 1$.  This implies that
$F$ is isomorphic to a subgroup of $H$.

\begin{cor}\label{cor:retract}
Under the situation of Proposition \emph{\ref{R:8.4.9}}, the free
group $F$ is a retract of $G(p)$. Hence if $J = \Ob(C)$ has more
than one element, then the group $G(p)$ is not trivial, and if $J$
has more than two elements, then $G(p)$ is not abelian.\end{cor}

\Proof Let $M = A(p) * B(p) * F$, and let $\imath' : F \to M$ be
the inclusion.  Let $\rho' : M \to F$ be the retraction which is
trivial on $A(p)$ and $B(p)$ and is the identity on $F$.  Let $q :
M \to G(p)$ be the quotient morphism.  Then it is clear that
$\rho'$ preserves the relations (\ref{Eq:8.4:2}), and so $\rho'$
defines uniquely a morphism $\rho : G(p) \to F$ such that $\rho q
= \rho'$.  Let $\imath = qi'$. Then $\rho \imath = \rho' i' = 1$.
So $F$ is a retract of $G(p)$.

The concluding statements are clear. \sqbox

We use the last two statements of the Corollary in sections 4 and
5 respectively.

\section{The Phragmen-Brouwer Property }

A topological space $X$ is said to \emph{have the Phragmen-Brouwer
Property} \index{Phragmen-Brouwer Property} (here abbreviated to
PBP) if $X$ is connected and the following holds: \emph{if $D$ and
$E$ are disjoint, closed subsets of $X$, and if $a$ and $b$ are
points in $X \setminus (D \cup E)$ which lie in the same component
of $X \setminus D$ and in the same component of $X \setminus E$,
then $a$ and $b$ lie in the same component of $X \setminus (D \cup
E)$.}  To express this more succinctly, we say a subset $D$ of a
space $X$ \emph{separates} \index{separates} the points $a$ and
$b$ if $a$ and $b$ lie in distinct components of $X \setminus D$.
Thus the PBP is that: if $D$ and $E$ are disjoint closed subsets
of $X$ and $a, b$ are points of $X$ not in $D \cup E$ such that
neither $D$ nor $E$ separate $a$ and $b$, then $D \cup E$ does not
separate $a$ and $b$.

A standard example of a space not having the PBP is the circle
$\ES^{1}$, since we can take $D = \{+1\}$, $E = \{-1\}$, $a = i$,
$b = -i$.  This example is typical, as the next result shows. But
first we remark that our criterion for the PBP will involve
fundamental groups, that is will involve paths, and so we need to
work with path-components rather than components. However, if $X$
is locally path-connected, then components and path-components of
open sets of $X$ coincide, and so for these spaces we can replace
in the PBP `component' by `path-component'. This explains the
assumption of locally path-connected in the results that follow.

\begin{prop} \label{R:8.5.1}  Let $X$  be a path-connected and locally
path-connected space whose fundamental group (at any point) does
not have the integers $\Z$ as a retract.  Then $X$ has the PBP.
\end{prop}

\Proof Suppose $X$ does not have the PBP.  Then there are
disjoint, closed subsets $D$ and $E$ of $X$ and points a and $b$
of $X\setminus(D \cup E)$ such that $D \cup E$ separates $a$ and
$b$ but neither $D$ nor $E$ separates $a$ and $b$.  Let $U = X
\setminus D$, $V = X \setminus E$, $W = X \setminus (D\cup E) = U
\cap V$.  Let $J$ be a subset of $W$ such that $a, b \in J$ and
$J$ meets each path-component of $W$ in exactly one point. Since
$D$ and $E$ do not separate $a$ and $b$, there are elements
$\alpha \in \pi_1  U(a,b)$ and $\beta \in \pi_1 V(a,b)$. Since $X$
is path-connected, the pairs $(U,J), (V,J), (W,J)$ are connected.
By the van Kampen Theorem \ref{thm:vKT} the following diagram of
morphisms induced by inclusions is a pushout of groupoids:
$$
\xymatrix@L=0pc@=1pc{ \pi_1  W J \ar[rr]^*+{i_1} \ar[dd]_*+{i_2}
&& \pi_1  U J \ar[dd]^*+{u_1} \\ \\ \pi_1  V J \ar[rr]_*+{u_2} &&
\pi_1  X J. }
$$

Since $U$ and $V$ are path-connected and $J$ has more than one
element, it follows from Corollary \ref{cor:retract} that $\pi_1
XJ$ has the integers $\Z$ as a retract. \sqbox

As an immediate application we obtain:

\begin{prop} \label{R:8.5.2}
 The following spaces have the PBP: the sphere $\ES^{n}$ for $n>1$;
$\ES^{2} \setminus \{a\}$ for $a \in \ES^{2}$; $\ES^{n} \setminus
\Lambda$ if $\Lambda$ is a finite set in $\ES^{n}$ and $n>2$.
\sqbox
\end{prop}

In each of these cases the fundamental group is trivial.

An important step in our proof of the Jordan Curve Theorem is to
show that if $A$ is an \emph{arc} in $\ES^{2}$, that is a subspace
of $\ES^{2}$ homeomorphic to the unit interval $\I$, then the
complement of $A$ is  path-connected.  This follows from the
following more general result.

\begin{prop} \label{R:8.5.3}
Let $X$ be a path-connected and locally path-connected Hausdorff
space such that for each $x$ in $X$ the space $X \setminus \{x\}$
has the PBP.  Then any arc in $X$ has path-connected complement.
\end{prop}

\Proof Suppose $A$ is an arc in $X$ and $X\setminus A$ is not
path-connected.  Let $a$ and $b$ lie in distinct path-components
of $X \setminus A$.

By choosing a homeomorphism $\I \to A$ we can speak unambiguously
of the mid-point of $A$ or of any subarc of $A$.  Let $x$ be the
mid-point of $A$, so that $A$ is the union of sub-arcs $A'$ and
$A''$ with intersection $\{x\}$.  Since $X$ is Hausdorff, the
compact sets $A'$ and $A''$ are closed in $X$.  Hence $A'
\setminus \{x\}$ and $A'' \setminus \{x\}$ are disjoint and closed
in $X \setminus \{x\}$.  Also $A\setminus \{x\}$ separates a and
$b$ in $X \setminus \{x\}$ and so one at least of $A' ,A''$
separates a and $b$ in $X \setminus \{x\}$.  Write $A_{1}$ for one
of $A', A''$ which does separate $a$ and $b$. Then $A_{1}$ is also
an arc in $X$.

In this way we can find by repeated bisection a sequence $A_{i}$,
$i\geqslant 1$, of sub-arcs of $A$ such that for all $i$ the
points $a$ and $b$ lie in distinct path-components of $X \setminus
A_{i}$ and such that the intersection of the $A_{i}$ for
$i\geqslant 1$ is a single point, say $y$, of $X$.

Now $X \setminus \{y\}$ is path-connected, by definition of the
PBP. Hence there is a path $\lambda$ joining $a$ to $b$ in $X
\setminus \{y\} $.  But $\lambda $ has compact image and hence
lies in some $X \setminus A_{i}$.  This is a contradiction. \sqbox

\begin{cor}\label{cor:compl-arc}
The complement of any arc in $\ES^{n}$ is path-connected.
\end{cor}
\noindent {\bf Sketch Proof} The case $n = 0$ is trivial, while the
case $n = 1$ needs a special argument that the complement of any arc
in $\ES^{1}$ is an open arc.  The case $n \geqslant 2$ follows from
the above results. \sqbox

\section{The Jordan Separation and Curve Theorems}

We now prove one step along the way to the full Jordan Curve
Theorem.

\begin{thm}[The Jordan Separation Theorem] \label{R:8.5.4}
\index{Jordan separation theorem} The complement of a simple
closed curve in $\ES^{2}$ is not connected.
\end{thm}

\Proof Let $C$ be a simple closed curve in $\ES^{2} $.  Since $C$
is compact and $\ES^{2}$ is Hausdorff, $C$ is closed, $\ES^{2}
\setminus C$ is open, and so path-connectedness of $\ES^{2}
\setminus C$ is equivalent to connectedness.

Write $C = A \cup B$ where $A$ and $B$ are arcs in $C$ meeting
only at $a$ and $b$ say.  Let $U = \ES^{2} \setminus A$, $V =
\ES^{2} \setminus B$, $W = U \cap V$, $X = U \cup V$.  Then $W =
\ES^{2} \setminus C$ and $X = \ES^{2} \setminus \{a,b\} $.  Also
$X$ is path-connected, and, by Corollary \ref{cor:compl-arc}, so
also are $U$ and $V$.

Let $x \in W$.  Suppose that $W$ is path-connected.  By the van
Kampen Theorem \ref{thm:vKT}, the following diagram of morphisms
induced by inclusion is a pushout of groups:
$$
\xymatrix@L=0pc@=1pc{ \pi_1  (W,x) \ar[rr] \ar[dd] &&
\pi_1  (U,x) \ar[dd]^*+{i_{*}} \\ \\
\pi_1  (V,x) \ar[rr]_*+{j_{*}} && \pi_1  (X,x). }
$$
Now $\pi_1  (X,x)$ is isomorphic to the group $\Z$ of integers.
We derive a contradiction by proving that the morphisms $i_{*}$
and $j_{*}$ are trivial.  We give the proof for $i_{*}$, as that
for $j_{*}$ is similar.

Let $f : \ES^{1} \to U$ be a map and let $g = if : \ES^{1} \to X$.
Let $\gamma$ be a parametrisation of $A$ which sends $0$ to $b$ and
$1$ to $a$. Choose a homeomorphism $h : \ES^{2} \setminus \{a\} \to
\Real^{2}$ which takes $b$ to $0$ and such that $hg$ maps $\ES^{1}$
into $\Real^{2} \setminus \{0\}$.  Then $h \gamma(0)=0$ and $\|h
\gamma(t)\| $ tends to infinity as $t$ tends to 1. Since the image
of $g$ is compact, there is an $r>0$ such that $hg[\ES^{1}]$ is
contained in $B(0,r)$. Now there exists $0 < t_0 < 1$ such that the
distance from $0$ to $y=h\gamma(t_0)$ is $>r$. Define the path
$\lambda$ to be the part of $h\gamma$ reparametrised so that
$\lambda(0)=0 $ and $ \lambda(1) =y$.

Define $G : \ES^{1} \times \I \to \Real^{2}$ by
$$
G(z,t) =
\begin{cases}
    hg(z)-\lambda (2t)  &\mbox{if}~0 \leqslant t \leqslant \frac{1}{2},\\
    (2-2t)hg(z)-y       &\mbox{if}~\frac{1}{2} \leqslant t \leqslant 1.
\end{cases}
$$
Then $G$ is well-defined.  Also $G$ never takes the value $0$
(this explains the choices of $\lambda$ and $y$).  So $G$ gives a
homotopy in $\Real^{2} \setminus \{0\}$ from $hg$ to the constant
map at $-y$. So $hg$ is inessential and hence $g$ is inessential.
This completes the proof that $i_{*}$ is trivial. \sqbox

As we shall see, the Jordan Separation Theorem is used in the
proof of the Jordan Curve Theorem.

\begin{thm} [Jordan Curve Theorem] \label{R:8.5.5}
\index{Jordan Curve Theorem} If $C$ is a simple closed curve in
$\ES^{2}$, then the complement of $C$ has exactly two components,
each with $C$ as boundary.
\end{thm}

\Proof As in the proof of Theorem \ref{R:8.5.4}, write $C$ as the
union of two arcs $A$ and $B$ meeting only at $a$ and $b$ say, and
let $U = \ES^{2} \setminus A$, $V = \ES^{2} \setminus B$.  Then
$U$ and $V$ are path-connected and $X = U \cup V = \ES^{2}
\setminus \{a,b\}$ has fundamental group isomorphic to $\Z$.  Also
$W = U \cap V = \ES^{2} \setminus C$ has at least two
path-components, by the Jordan Separation Theorem  \ref{R:8.5.4}.

If $W$ has more than two path-components, then the fundamental
group $G$ of $X$ contains a copy of the free group on two
generators, by Corollary \ref{cor:retract}, and so $G$ is
non-abelian. This is a contradiction, since $G \cong \Z$.  So $W$
has exactly two path-components $P$ and $Q$, say, and this proves
the first part of Theorem \ref{R:8.5.5}.

Since $C$ is closed in $\ES^{2}$ and $\ES^{2}$ is locally
path-connected, the sets $P$ and $Q$ are open in $\ES^{2}$.  It
follows that if $x \in \overline{P} \setminus P$ then $x \notin
Q$, and hence $\overline{P} \setminus P$ is contained in $C$.  So
also is $\overline{Q} \setminus Q$, for similar reasons.  We prove
these sets are equal to $C$.

Let $x \in C$ and let $N$ be a neighbourhood of $x$ in $\ES^{2}$.
We prove $N$ meets $\overline{P} \setminus P$.  Since
$\overline{P} \setminus P$ is closed and $N$ is arbitrary, this
proves that $x \in \overline{P} \setminus P$.

Write $C$ in a possibly new way as a union of two arcs $D$ and $E$
intersecting in precisely two points and such that $D$ is contained
in $N \cap C$.  Choose points $p$ in $P$ and $q$ in $Q$.  Since
$\ES^{2} \setminus E$ is path-connected, there is a path $\lambda$
joining $p$ to $q$ in $\ES^{2} \setminus E$.  Then $\lambda$ must
meet $D$, since $p$ and $q$ lie in distinct path-components of
$\ES^{2} \setminus E$.  In fact if $s = \sup\{t \in \I : \lambda
[0,t]\subseteq P\} $, then $\lambda (s) \in \overline{P} \setminus
P$.  It follows that $N$ meets $\overline{P} \setminus P$.

So $\overline{P} \setminus P = C$ and similarly $\overline{Q}
\setminus Q = C$. \sqbox \providecommand{\BIBtom}{tom}
\providecommand{\BIBvan}{Van}

\begin{small}
  \printindex
\end{small}

\end{document}